\documentclass[12pt]{amsart}

\usepackage{amssymb}

\newtheorem{theorem}{Theorem}[section]
\newtheorem{proposition}[theorem]{Proposition}
\newtheorem{lemma}[theorem]{Lemma}
\newtheorem{corollary}[theorem]{Corollary}
\newtheorem{definition}[theorem]{Definition}
\newtheorem{remark}[theorem]{Remark}

\numberwithin{equation}{section}

\begin{document}

\baselineskip=15.5pt

\title[Semistable bundles over projective manifolds, II]{On
semistable principal bundles over a complex projective 
manifold, II}

\author[I. Biswas]{Indranil Biswas}

\address{School of Mathematics, Tata Institute of Fundamental
Research, Homi Bhabha Road, Bombay 400005, India}

\email{indranil@math.tifr.res.in}

\author[U. Bruzzo]{Ugo Bruzzo}

\address{Scuola Internazionale Superiore di Studi Avanzati,
Via Beirut 2--4, 34013, Trieste, Italy, and Istituto Nazionale di 
Fisica Nucleare, Sezione di Trieste}

\email{bruzzo@sissa.it}



\date{}

\begin{abstract}
Let $(X,\,\omega)$ be a compact
connected K\"ahler manifold of complex
dimension $d$ and $E_G\,\longrightarrow\, X$ a holomorphic 
principal $G$--bundle, where $G$
is a connected reductive linear algebraic group
defined over $\mathbb C$. Let $Z(G)$ denote
the center of $G$. We prove that the
following three statements are equivalent:
\begin{enumerate}
\item{} There is a parabolic subgroup $P\, \subset\, G$
and a holomorphic reduction of structure group
$E_P\, \subset\, E_G$ to $P$, such that the corresponding
$L(P)/Z(G)$--bundle
$$
E_{L(P)/Z(G)}\, :=\, E_P(L(P)/Z(G))\, \longrightarrow\, X
$$
admits a unitary flat connection, where
$L(P)$ is the Levi quotient of $P$.

\item{} The adjoint vector bundle ${\rm ad}(E_G)$ is
numerically flat.

\item{} The principal $G$--bundle $E_G$ is pseudostable,
and
$$
\int_X c_2({\rm ad}(E_G))\omega^{d-2}\, =\, 0\, .
$$
\end{enumerate}

If $X$ is a complex projective manifold, and $\omega$
represents a rational cohomology class, then the third
statement is equivalent to the statement that
$E_G$ is semistable with $c_2({\rm ad}(E_G))\, =\, 0$.
\end{abstract}

\maketitle

\section{Introduction}

Let $G$ be a connected reductive linear algebraic group defined
over $\mathbb C$. Let $Z(G)\, \subset\, G$ be the center of $G$.
Fix a proper parabolic subgroup
$Q\, \subset\, G$ without any simple factor.
Also, fix a character $\chi$ of $Q$ that satisfies the following
two conditions:
\begin{itemize}
\item{} $\chi$ is trivial on
$Z(G)$, and

\item{} for each simple factor $H$ of $G/Z(G)$,
the restriction of $\chi$ to the parabolic subgroup
$H\bigcap (Q/Z(G))\, \subset\, H$
is nontrivial and antidominant.
\end{itemize}

Let $E_G$ be a holomorphic principal $G$--bundle over a connected
complex projective manifold $M$. Fix a polarization on $M$.
We have the line bundle
$$
L_\chi \, :=\, (E_G\times {\mathbb C})/Q\,
\longrightarrow\, E_G/Q
$$
associated to the principal $Q$--bundle
$E_G\, \longrightarrow\, E_G/Q$ for the character $\chi$.

The following theorem was proved in \cite{BB} (see
\cite[p. 24, Theorem 4.3]{BB}):

\begin{theorem}\label{thm-gb}
The following four statements are equivalent.
\begin{enumerate}
\item{} The principal $G$--bundle $E_G$ is semistable, and the
second Chern class
$$
c_2({\rm ad}(E_G)) \, \in\, H^4(M,\, {\mathbb Q})
$$
vanishes.

\item{} The associated line bundle $L_\chi$ over $E_G/Q$
is numerically effective.

\item{} For every pair of the form $(Y\, ,\psi)$,
where $Y$ is a compact connected Riemann surface and
$$
\psi\,:\, Y\,\longrightarrow\, M
$$
a holomorphic map, and for every
holomorphic reduction $E_Q\, \subset\,
\psi^*E_G$ of structure group to $Q$ of the principal
$G$--bundle $\psi^*E_G\,\longrightarrow\,Y$, the associated
line bundle
$$
E_Q(\chi)\, =\, (E_Q\times {\mathbb C})/Q\,\longrightarrow\,Y
$$
is of nonnegative degree.

\item{} For any pair $(Y\, ,\psi)$ as in statement (3), the
principal $G$--bundle $\psi^*E_G$ over $Y$ is semistable.
\end{enumerate}
\end{theorem}

We also recall that a vector bundle $E$ over a compact K\"ahler
manifold is polystable with $c_i(E)\,=\, 0$ for all
$i\,\geq\, 1$ if and only if
$E$ admits a unitary flat connection \cite{Do}, \cite{UY}. This
holds also for principal bundles \cite{AB}.

We show that the equivalent conditions in Theorem \ref{thm-gb}
can be expressed using unitary flat connections
on some associated bundles.

For a parabolic subgroup $P\, \subset\, G$,
its Levi quotient will be denoted by $L(P)$. The center $Z(G)$
is contained in $P$, and it projects isomorphically to $L(P)$.
So $Z(G)$ will be considered as a subgroup of $L(P)$.

We prove the following theorem (see Theorem \ref{thm1}):

\begin{theorem}\label{thm0}
Let $E_G$ be a holomorphic principal $G$--bundle over
a compact connected K\"ahler manifold $(X\, ,\omega)$, where
$G$ is a connected reductive linear algebraic group defined
over $\mathbb C$. Then the following three statements are
equivalent:
\begin{enumerate}

\item{} There is a parabolic subgroup $P\, \subset\, G$
and a holomorphic reduction of structure group
$E_P\, \subset\, E_G$ to $P$, such that the corresponding
$L(P)/Z(G)$--bundle
$$
E_P(L(P)/Z(G))\,=\, (E_P\times (L(P)/Z(G)))/P
\, \longrightarrow\, X
$$
admits a unitary flat connection (see Definition
\ref{def.unit.con}).

\item{} The adjoint vector bundle ${\rm ad}(E_G)$ is
numerically flat.

\item{} The principal $G$--bundle $E_G$ is pseudostable
(see Definition \ref{d-p}),
and
$$
\int_X c_2({\rm ad}(E_G))\omega^{d-2}\, =\, 0\, .
$$
\end{enumerate}
\end{theorem}

If $X$ is a complex projective manifold, and $\omega$
represents a rational cohomology class, then the third
statement in Theorem \ref{thm1} coincides with the
first statement in Theorem \ref{thm-gb} (see Corollary
\ref{cor2}). An equivalent form of this result was proved in 
\cite{BG2} (as a particular
case of a result valid for principal Higgs bundles).

\section{Pseudostability and numerically flatness}\label{sec2}

Let $X$ be a compact connected K\"ahler manifold
of complex dimension $d$. Fix a
K\"ahler form $\omega$ on $X$.
Let $d$ be the complex dimension of $X$.
For any torsionfree coherent analytic sheaf $F$ on $X$, the
\textit{degree} of $F$ is defined to be
$$
\text{degree}(F)\,:=\, \int_X c_1(F)\wedge
\omega^{d-1}\, \in\, {\mathbb R}\, .
$$
Let $U\, \subset\, X$ be any dense open subset
such that the complement $X\setminus U$ is a complex analytic
subspace of $X$ of complex codimension at least two.
Let $\iota\, :\, U\, \hookrightarrow\, X$ be the inclusion
map. For any holomorphic vector bundle $E$ over $U$ such that
the direct image $\iota_* E$ is a coherent analytic sheaf
on $X$, define
$$
\text{degree}(E)\, :=\, \text{degree}(\iota_* E)\, .
$$

We note that
in the special case where $X$ is a complex projective manifold,
for any vector bundle $E$ defined over a Zariski open dense
subset $U\, \subset\, X$ such that complement $X\setminus U$ is
of complex codimension at least two, the direct image
$\iota_* E$ is a coherent sheaf on $X$, where $\iota$ as before
is the inclusion of $U$ in $X$.

If $E\, \not=\, 0$, it is customary to denote the real
number $\text{degree}(E)/\text{rank}(E)$ by $\mu(E)$; this
number $\mu(E)$ is also called the \textit{slope} of $E$.

A torsionfree coherent analytic sheaf
$F$ over $X$ is called \textit{stable} (respectively,
\textit{semistable}) if
$$
\mu(F')\, <\, \mu(F)
$$
(respectively, $\mu(F') \, \leq\, \mu(F)$) for all
coherent subsheaves $F'\,\subset\, F$
with $0\, <\, \text{rank}(F')\, <\, \text{rank}(F)$.
A semistable sheaf is called \textit{polystable} if
it is a direct sum of stable sheaves.

A holomorphic vector bundle $F$ over $X$ is called
\textit{pseudostable} if there exists a filtration of
holomorphic subbundles
of $F$
\begin{equation}\label{fil1}
0\, =\, F_0\, \subset\, F_1 \, \subset\, \cdots
\, \subset\, F_{n-1} \, \subset\, F_n \, =\, F
\end{equation}
such that for each $i\, \in\, [1\, ,n]$, the quotient
$F_i/F_{i-1}$ is a polystable vector bundle with
$\mu(F_i/F_{i-1}) \, =\,\mu(F)$.

Therefore, any pseudostable vector bundle is semistable.
We note that any semistable vector bundle $E$ admits a
filtration of coherent analytic subsheaves
$$
0\, =\, E_0\, \subset\, E_1 \, \subset\, \cdots
\, \subset\, E_{m-1} \, \subset\, E_m \, =\, E
$$
such that the quotient sheaf $E_i/E_{i-1}$ is polystable
with $\mu(E_i/E_{i-1}) \, =\,\mu(E)$ for all
$i\, \in\, [1\, ,m]$ (see \cite[p. 23, Lemma 1.5.5]{HL}).
However, the subsheaves
$F_i$ in Eq. \eqref{fil1} are required to be locally free.

We will briefly recall the definitions of stable
and semistable principal bundles
(see \cite{Ra}, \cite{AB} for the details).

Let $G$ be a connected reductive linear algebraic group
defined over the field of complex numbers. Let $E_G$ be a
holomorphic principal $G$--bundle over $X$.

Consider
all triples $(P\, ,U\, , \sigma)$ of the following form:
\begin{itemize}
\item $P\, \subset\, G$ is a proper maximal parabolic subgroup,

\item $\iota\, :\, U\, \hookrightarrow\, X$ is a dense open
subset such that the complement $X\setminus U$
is an analytic subset of $X$ of (complex) codimension at
least two, and

\item $\sigma\, :\, U\, \longrightarrow\,(E_G/P)\vert_U$
is a holomorphic
reduction of structure group of $E_G$ to $P$ over $U$
satisfying the condition that the direct image
$\iota_* \sigma^* T_{\text{rel}}$ is a coherent analytic
sheaf on $X$, where $T_{\text{rel}}\, \longrightarrow\,
E_G/P$ is the relative 
tangent
bundle for the projection $E_G/P\, \longrightarrow\, X$.
\end{itemize}

The principal
$G$--bundle $E_G$ is called \textit{semistable} (respectively,
\textit{stable}) if for each
triple $(P\, ,U\, , \sigma)$ of the above type, the inequality
$$
\text{degree}(\sigma^* T_{\text{rel}}) \, \geq \, 0
$$
(respectively, $\text{degree}(\sigma^* T_{\text{rel}}) \,
> \, 0$) holds.

Let
\begin{equation}\label{zg}
Z(G)\, \subset\, G
\end{equation}
be the center. We note that $Z(G)$ is contained in all
parabolic subgroups of $G$.

A holomorphic reduction of structure group of $E_G$
$$
E_Q\, \subset\, E_G
$$
to some parabolic subgroup $Q\, \subset\, G$ is called
\textit{admissible} if for each character $\lambda$ of
$Q$ which is trivial on $Z(G)$, the degree of
the line bundle $E_Q({\mathbb C})\,=\, (E_Q\times
{\mathbb C})/Q\,\longrightarrow\, X$ associated to $Q$ for
$\lambda$ is zero \cite[p. 307, Definition 3.3]{Ra}.

A holomorphic principal $G$--bundle $E_G\,\longrightarrow\, X$
is called \textit{polystable} if either $E_G$ is stable, or there
is a proper parabolic subgroup $Q\, \subset\, G$
and a holomorphic reduction of structure group
$$
E_G \, \supset \, E_{L(Q)}\, \longrightarrow\, X
$$
to a Levi subgroup $L(Q)$ of $Q$ such that the
following two conditions hold:
\begin{enumerate}
\item the principal $L(Q)$--bundle $E_{L(Q)}$ is stable, and

\item the reduction of structure group of $E_G$ to $Q$
obtained by extending the structure group of $E_{L(Q)}$
using the inclusion of $L(Q)$ in $Q$ is admissible.
\end{enumerate}
(A Levi subgroup of a parabolic group $Q$ is a
maximal reductive subgroup of $Q$. Any Levi subgroup
of $Q$ projects isomorphically onto the Levi quotient
of $Q$; see \cite[p. 217, Theorem 7.1]{Mo}.)

For a principal $G$--bundle $E_G\, \longrightarrow\,
X$, its adjoint vector bundle
$$
\text{ad}(E_G)\, :=\, E_G({\mathfrak g})\, =\,
(E_G\times {\mathfrak g})/G \, \longrightarrow\, X
$$
is the one associated to $E_G$ for the adjoint action of $G$
on its Lie algebra $\mathfrak g$. We recall from \cite{BG}
the definition of a pseudostable principal $G$--bundle.

\begin{definition}{\rm (}\cite[p. 26,
Definition 2.3]{BG}{\rm ).}\label{d-p}
{\rm A holomorphic principal $G$--bundle $$E_G\,\longrightarrow
\, X$$ is called} pseudostable {\rm if its
adjoint vector bundle $\text{ad}(E_G)$ is pseudostable.}
\end{definition}

A principal $G$--bundle $E_G$ is semistable if and only
if the adjoint vector bundle $\text{ad}(E_G)$ is semistable
\cite[p. 214, Proposition 2.10]{AB}.
Therefore, a pseudostable principal bundle is semistable.

A holomorphic
line bundle $L$ over a compact complex manifold $M$
is called \textit{numerically effective} if
$L$ admits Hermitian structures
whose curvatures have arbitrarily small negative part
(see \cite[p. 299, Definition 1.2]{DPS}). When $M$ is a
complex projective manifold, this definition coincides with the
usual definition of numerically effectiveness which says that
the degree of the restriction of the line bundle to each
irreducible complete curve is nonnegative.

A holomorphic vector bundle $E$ over $M$ is called
\textit{numerically effective} if the tautological
line bundle ${\mathcal O}_{{\mathbb P}(E)}(1)$ over the total
space of the projective bundle ${\mathbb P}(E)\,\longrightarrow
\, M$ is numerically effective \cite[p. 305, Definition 1.9]{DPS}.
A holomorphic vector bundle $E$ over $M$ is called
\textit{numerically flat} if both $E$ and $E^*$ are
numerically effective \cite[p. 311, Definition 1.17]{DPS}.

The following theorem is proved in \cite{DPS}
(see \cite[p. 311, Theorem 1.18]{DPS}):

\begin{theorem}[\cite{DPS}]\label{thm-dps}
A holomorphic vector bundle $E\, \longrightarrow\, X$
is numerically flat if and only if $E$ admits a filtration
of holomorphic subbundles
such that each successive quotient admits a unitary
flat connection.
\end{theorem}

\begin{lemma}\label{lem0}
Let $E\, \longrightarrow\, X$ be a numerically flat vector
bundle. Then $E$ is pseudostable, and $c_i(E)\, \in\,
H^{2i}(X,\, {\mathbb Q})$ vanishes for all $i\, \geq\, 1$.
\end{lemma}

\begin{proof}
A holomorphic vector bundle $V\, \longrightarrow\, X$
admits a unitary flat connection if and only if $V$
is polystable and $c_i(V)\, =\, 0$ for all $i\, \geq\, 1$
\cite{Do}, \cite{UY}.
Therefore, the lemma follows from Theorem \ref{thm-dps}.
\end{proof}

Let $H$ be a reductive linear algebraic group defined over
$\mathbb C$. Fix a maximal compact subgroup $K(H)\,\subset
\, H$.

\begin{definition}\label{def.unit.con}
{\rm A holomorphic principal $H$--bundle $E_H\,\longrightarrow
\, X$ is said to admit a} unitary flat connection {\rm if
$E_H$ is given by a representation of the fundamental group
of $X$ in $K(H)$.}
\end{definition}

\begin{remark}\label{adj.}
{\rm Since $H$ is reductive, the Lie algebra of
$H$, considered as a $H$--module, is self--dual.
Therefore, for any principal $H$--bundle $E_H$, we
have ${\rm ad}(E_H)\,=\, {\rm ad}(E_H)^*$.}
\end{remark}

\section{Projectively unitary flat Levi bundles}

Let $G$ be a connected reductive linear algebraic group
defined over $\mathbb C$. The Lie algebra
of $G$ will be denoted by $\mathfrak g$. Let
$$
P\, \subset\, G
$$
be a parabolic subgroup. Let
\begin{equation}\label{rup}
R_u(P)\, \subset\, P
\end{equation}
be the unipotent radical. So $R_u(P)$ is the unique maximal
connected normal unipotent subgroup of $P$.
Consider the adjoint action of $P$ on the Lie algebra
$\mathfrak g$. Using it, $\mathfrak g$ will be considered
as a $P$--module.

\begin{lemma}\label{lem1}
There is a filtration of $P$--modules
$$
0\, =\, W_0\, \subset\, W_1 \, \subset\, \cdots \, \subset\,
W_{\ell-1} \, \subset\, W_\ell\, =:\, \mathfrak g
$$
such that for each $i\, \in\, [1\, ,\ell]$, the action on
$W_i/W_{i-1}$ of the subgroup $R_u(P)$ is the trivial one.
\end{lemma}

\begin{proof}
Set
$$
W_1\, :=\, \{v\,\in\, {\mathfrak g}\, \mid\, \text{Ad}(g)(v)
\, =\, v ~\ \forall~\ g\,\in\, R_u(P)\}\, .
$$
Since $R_u(P)$ is a normal subgroup of $P$, the action of $P$
on $\mathfrak g$ preserves $W_1$. Now construct $W_i$ inductively
in the following way.

Define
$$
W'_j\, :=\, \{v\,\in\, {\mathfrak g}/W_{j-1}\, \mid\,
\text{Ad}(g)(v) \, =\, v ~\ \forall~\ g\,\in\, R_u(P)\}\, .
$$
The action of $P$ on ${\mathfrak g}/W_{j-1}$ preserves $W'_j$.
Now set
$$
W_j\, :=\, q^{-1}_{j-1}(W'_j)\, \subset\, {\mathfrak g}\, ,
$$
where $q_{j-1}\, :\, {\mathfrak g}\,\longrightarrow\,{\mathfrak 
g}/W_{j-1}$ is the quotient map. The resulting filtration
$$
0\, =\, W_0\, \subset\, W_1 \, \subset\, \cdots \, \subset\,
W_{\ell-1} \, \subset\, W_\ell\, =:\, \mathfrak g
$$
clearly satisfies the condition in the lemma.
\end{proof}

Let $E_G\, \longrightarrow\, X$ be a holomorphic principal
$G$--bundle.

Let $P\, \subset\, G$ be a parabolic subgroup. Let
$$
L(P)\, :=\, P/R_u(P)
$$
be the Levi quotient, which is a connected reductive group.
Consider the quotient map $P\, \longrightarrow\, L(P)$. Its
restriction to the center $Z(G)$ (see Eq. \eqref{zg}) is an
isomorphism. Using this map, we will consider $Z(G)$ as
a subgroup of $L(P)$. Note that $Z(G)$ 
is contained in the center of $L(P)$. Let
\begin{equation}\label{ep0}
L'(P)\, :=\, L(P)/Z(G)
\end{equation}
be the quotient group. For a principal $P$--bundle
$E_P\, \longrightarrow\, X$, let
\begin{equation}\label{ep-1}
E_{L'(P)}\, :=\, E_P(L'(P))\,=\, (E_P\times L'(P))/P
\, \longrightarrow\, X
\end{equation}
be the principal $L'(P)$--bundle obtained by extending the
structure group of $E_P$ using the quotient map of $P$ to
$L'(P)$ defined in Eq. \eqref{ep0}.

\begin{proposition}\label{prop1}
Assume that there is a holomorphic reduction of structure group
\begin{equation}\label{ep}
E_P\, \subset\, E_G
\end{equation}
satisfying the following conditions:
\begin{enumerate}
\item the principal $L'(P)$--bundle $E_{L'(P)}$
defined in Eq. \eqref{ep-1} is polystable,

\item for every character $\chi$ of $L'(P)$, the line
bundle
\begin{equation}\label{epchi}
E_{L'(P)}(\chi) \,:=\, (E_{L'(P)}\times
{\mathbb C})/L'(P)\, \longrightarrow\, X
\end{equation}
associated to $E_{L'(P)}$ for $\chi$ is of degree zero, and
\begin{equation}\label{epc2}
\int_X c_1(E_{L'(P)}(\chi))^2\omega^{d-2}\, =\, 0\, ,
\end{equation}
where $d\,=\, \dim_{\mathbb C} X$, and

\item for the adjoint vector bundle ${\rm ad}(E_{L'(P)})$
of $E_{L'(P)}$,
\begin{equation}\label{c2}
\int_X c_2({\rm ad}(E_{L'(P)}))\omega^{d-2}\, =\, 0\, .
\end{equation}
\end{enumerate}

Then the adjoint vector bundle ${\rm ad}(E_G)$ admits
a filtration of subbundles
$$
0\, =\, V_0\, \subset\, V_1 \, \subset\, \cdots \, \subset\,
V_{\ell-1} \, \subset\, V_\ell\, =\, {\rm ad}(E_G)
$$
such that for each $i\, \in\, [1\, ,\ell]$, the quotient
vector bundle $V_i/V_{i-1}$ admits a unitary flat connection.
\end{proposition}

\begin{proof}
Consider the filtration $\{W_i\}^\ell_{i=0}$ of
$\mathfrak g$ in Lemma \ref{lem1}. Let
$$
E_{W_i}\, :=\, E_P(W_i)\, \longrightarrow\, X
$$
be the holomorphic vector bundle associated to the
principal $P$--bundle $E_P$ in Eq. \eqref{ep} for the
$P$--module $W_i$. So the filtration $\{W_i\}^\ell_{i=0}$
gives a filtration of holomorphic vector bundles
\begin{equation}\label{fih}
0\, =\, E_{W_0}\, \subset\, E_{W_1} \, \subset\, \cdots \, 
\subset\,
E_{W_{\ell-1}} \, \subset\, E_{W_\ell}\, =\,
\text{ad}(E_G)\, .
\end{equation}
For any $i\,\in\, [1\, ,\ell]$,
the quotient vector bundle $E_{W_i}/E_{W_{i-1}}$ is clearly
identified with the vector bundle associated to the
principal $P$--bundle $E_P$ for the
$P$--module $W_i/W_{i-1}$.

Take any $i\,\in\, [1\, ,\ell]$. Since $R_u(P)$
acts trivially on $W_i/W_{i-1}$, the action of
$P$ on $W_i/W_{i-1}$ factors through an action of
the Levi quotient $L(P)\, =\, P/R_u(P)$
on $W_i/W_{i-1}$. Also, the
subgroup $Z(G)$ in Eq. \eqref{zg} acts trivially
on $\mathfrak g$. Consequently, the action of
$P$ on $W_i/W_{i-1}$ factors through an action of the quotient
$L'(P)$ in Eq. \eqref{ep0}. The quotient bundle
$E_{W_i}/E_{W_{i-1}}$ is identified with the vector
bundle $E_{L'(P)}(W_i/W_{i-1})$
associated to the principal $L'(P)$--bundle
$E_{L'(P)}$ in Eq. \eqref{ep-1} for the
$L'(P)$--module $W_i/W_{i-1}$.

Since $E_{L'(P)}$ is polystable, it admits a unique
Hermitian--Einstein connection \cite[p. 208, Theorem 0.1]{AB};
we will denote this connection by $\nabla$.
Consider the Hermitian--Einstein connection on the adjoint
vector bundle $\text{ad}(E_{L'(P)})$ induced by the
connection $\nabla$ on $E_{L'(P)}$. Since
$\bigwedge^{\text{top}} \text{ad}(E_{L'(P)})$ is a trivial
line bundle (see Remark \ref{adj.}),
and Eq. \eqref{c2} holds, it follows that the 
Hermitian--Einstein connection
on $\text{ad}(E_{L'(P)})$ is flat
\cite[Ch.~IV, \S~4, p. 114, Theorem (4.7)]{Ko}.

Given any character $\chi$ of $L'(P)$, consider the
Hermitian--Einstein connection on the associated
line bundle $E_{L'(P)}(\chi)$ (see Eq. \eqref{epchi})
induced by $\nabla$. Since $\text{degree}(E_{L'(P)}(\chi))
\,=\, 0$ (see condition (2) in
the proposition), and Eq. \eqref{epc2} holds, from
\cite[Ch.~IV, \S~4, p. 114, Theorem (4.7)]{Ko} it
follows that the Hermitian--Einstein connection on
$E_{L'(P)}(\chi)$ is flat. We noted above that
the induced connection on $\text{ad}(E_{L'(P)})$ is flat.
The quotient $L'(P)/[L'(P)\, ,L'(P)]$ is a product of
copies of ${\mathbb G}_m$, and the homomorphism
$$
L'(P)\, \longrightarrow\, \text{GL}(\text{Lie}(L'(P)))\times
(L'(P)/[L'(P)\, ,L'(P)])
$$
gives an injective homomorphism of Lie algebras;
here $L'(P)\, \longrightarrow\, \text{GL}(\text{Lie}(L'(P)))$
is the adjoint action.
Therefore, from the above observations that the
induced connections on $\text{ad}(E_{L'(P)})$ and
all associated line bundles are flat we conclude that
the Hermitian--Einstein connection $\nabla$ on
$E_{L'(P)}$ is flat.

Consequently, the connection on the associated vector
bundle $E_{L'(P)}(W_i/W_{i-1})$ induced by $\nabla$ is 
unitary flat. In view of the filtration in Eq.
\eqref{fih}, this completes the proof of the proposition.
\end{proof}

Combining Theorem \ref{thm-dps} and Proposition \ref{prop1},
we have the following corollary:

\begin{corollary}\label{cor1}
Assume that $E_G$ satisfies the conditions in Proposition
\ref{prop1}. Then the vector bundle ${\rm ad}(E_G)$ is
numerically flat.
\end{corollary}

\section{Pseudostable adjoint bundle}

Let
\begin{equation}\label{e0}
E_G\, \longrightarrow\, X
\end{equation}
be a pseudostable principal $G$--bundle. Define
\begin{equation}\label{e0p}
G'\, :=\, G/Z(G)\, ,
\end{equation}
where $Z(G)$ is the center (see Eq. \eqref{zg}). Let
\begin{equation}\label{e00}
E_{G'}\, :=\, E_G(G')\, =\, E_G/Z(G)
\, \longrightarrow\, X
\end{equation}
be the principal $G'$--bundle obtained by extending the
structure group of the principal $G$--bundle in Eq. \eqref{e0}
using the quotient map $G\, \longrightarrow\, G'$. Let
$\text{ad}(E_{G'})$ be the adjoint vector bundle of $E_{G'}$.
The Lie algebra of $G'$ will be denoted by ${\mathfrak g}'$.

Note that we have a decomposition of $G$--modules
$$
{\mathfrak g}\, =\, {\mathfrak g}'\oplus {\mathfrak z}
({\mathfrak g})\, ,
$$
where ${\mathfrak z}
({\mathfrak g})$ is the Lie algebra of $Z(G)$. Therefore,
\begin{equation}\label{d}
\text{ad}(E_G)\, =\, \text{ad}(E_{G'})\oplus (X\times
{\mathfrak z} ({\mathfrak g}))\, ,
\end{equation}
where $X\times {\mathfrak z} ({\mathfrak g})$ is the trivial
vector bundle over $X$ with fiber ${\mathfrak z}
({\mathfrak g})$. Since $\text{ad}(E_G)$ is
pseudostable, from Eq. \eqref{d} we conclude that the vector 
bundle $\text{ad}(E_{G'})$ is also pseudostable.
In particular, $\text{ad}(E_{G'})$ is semistable.

Note that $\text{degree}(\text{ad}(E_{G'}))
\,=\, 0$ (see Remark \ref{adj.}). Let
$$
S_1\, \subset\, \text{ad}(E_{G'})
$$
be the unique maximal polystable subsheaf of degree zero
of the semistable vector bundle $\text{ad}(E_{G'})$; this
subsheaf $S_1$ is called the \textit{socle} of $\text{ad}(E_{G'})$
(see \cite[p. 23, Lemma 1.5.5]{HL}). The condition
that $\text{ad}(E_{G'})$ is pseudostable implies that
$S_1$ is in fact a subbundle of $\text{ad}(E_{G'})$.
We note that the quotient $\text{ad}(E_{G'})/S_1$
is also a pseudostable vector bundle of degree zero.

Now inductively define $S'_{i+1}$ to be the socle of 
$\text{ad}(E_{G'})/S_i$. Inductively, it follows that
$S'_{i+1}$ is locally free. Now define
$$
S_{i+1}\, =\, \eta^{-1}_i(S'_{i+1})\, \subset\,
\text{ad}(E_{G'})\, ,
$$
where $\eta_i\, :\, \text{ad}(E_{G'})\, \longrightarrow\,
\text{ad}(E_{G'})/S_i$ is the quotient map. Therefore, we have
a filtration of subbundles
\begin{equation}\label{S}
0\, =\, S_0\, \subset\, S_1\, \subset\, \cdots\, \subset\,
S_{\ell-1} \, \subset\, S_\ell \, =\, \text{ad}(E_{G'})
\end{equation}
such that each successive quotient is a polystable vector
bundle of degree zero.

The filtration in Eq. \eqref{S} gives a holomorphic reduction
of structure group
\begin{equation}\label{equ}
E_Q\, \subset\, E_{G'}
\end{equation}
over $X$, where $Q$ is a parabolic subgroup of $G'$; see
\cite[p. 218]{AB}. We note that the reduction of structure
group in \cite[p. 218]{AB}
is only over an open subset of $X$ because the filtration of
the adjoint vector bundle in \cite{AB}
is by subsheaves and not necessarily by
subbundles; while in Eq. \eqref{S}, each $S_i$ is a subbundle
of $\text{ad}(E_{G'})$. Hence the reduction of structure group
$E_Q$ in Eq. \eqref{equ} is over entire $X$. The integer
$\ell$ in Eq. \eqref{S} is odd. The reduction $E_Q$ is
defined by the following condition: the adjoint vector bundle
$$
\text{ad}(E_Q)\, \subset\, \text{ad}(E_{G'})
$$
coincides with $S_{(\ell+1)/2}$. This condition determines
the pair $(Q\, , E_Q)$ in the following sense:
for any other pair $(Q_1\, , E_{Q_1})$ satisfying
this condition, there is some $g\, \in\, G'$ such that
\begin{itemize}
\item $Q_1\, =\, g^{-1}Qg$, and

\item $E_{Q_1}\, =\, E_Qg$.
\end{itemize}
See \cite{AB} for the details.

Fix a pair $(Q\, ,E_Q)$ such that $\text{ad}(E_Q)\, \subset\, 
\text{ad}(E_{G'})$ coincides with $S_{(\ell+1)/2}$.

Let
\begin{equation}\label{Q}
R_u(Q)\, \subset\, Q
\end{equation}
be the unipotent radical of $Q$. Let
\begin{equation}\label{LQ}
L(Q)\, :=\, Q/R_u(Q)
\end{equation}
be the Levi quotient.

Let
\begin{equation}\label{EL}
E_{L(Q)}\, :=\, E_Q(L(Q)) \, \longrightarrow\, X
\end{equation}
be the principal $L(Q)$--bundle obtained by extending the
structure group of $E_Q$ (see Eq. \eqref{equ}) using the
quotient map of $Q$ to $L(Q)$ in Eq. \eqref{LQ}.

\begin{proposition}\label{prop2}
Let $E_G\, \longrightarrow\, X$ be a pseudostable principal
$G$--bundle such that
\begin{equation}\label{a}
\int_X c_2({\rm ad}(E_{G'}))\omega^{d-2}\, =\, 0\, ,
\end{equation}
where $d\,=\, \dim_{\mathbb C} X$ and $E_{G'}$ is
defined in Eq. \eqref{e00}.

Then the principal $L(Q)$--bundle $E_{L(Q)}$ in
Eq. \eqref{EL} is polystable.

Also, the following two statements hold:
\begin{enumerate}
\item For every character $\chi$ of $L(Q)$, the line
bundle
$$
E_{L(Q)}(\chi) \, \longrightarrow\, X
$$
associated to $E_{L(Q)}$ for $\chi$ is of degree zero, and
\begin{equation}\label{c03}
\int_X c_1(E_{L(Q)}(\chi))^2\omega^{d-2}\, =\, 0\, .
\end{equation}

\item For the adjoint vector bundle ${\rm ad}(E_{L(Q)})$
of $E_{L(Q)}$,
\begin{equation}\label{c3}
\int_X c_2({\rm ad}(E_{L(Q)}))\omega^{d-2}\, =\, 0\, .
\end{equation}
\end{enumerate}
\end{proposition}

\begin{proof}
As we mentioned earlier, the subbundle $\text{ad}(E_Q)\,
\subset\, \text{ad}(E_{G'})$ coincides with $S_{(\ell+1)/2}$
in Eq. \eqref{S}. The unipotent radical subbundle of
$\text{ad}(E_Q)$ coincides with $S_{(\ell-1)/2}$. Hence
the adjoint vector bundle $\text{ad}(E_{L(Q)})$ is identified
with the quotient bundle $S_{(\ell+1)/2}/S_{(\ell-1)/2}$
(see \cite{AB}). Since
each successive quotient in Eq. \eqref{S} is polystable, we
conclude that $\text{ad}(E_{L(Q)})$ is a polystable
vector bundle. Therefore,
the principal $L(Q)$--bundle $E_{L(Q)}$ is polystable
\cite[p. 224, Corollary 3.8]{AB}.

We will now show that Eq. \eqref{c3} holds.

Let $W\, \longrightarrow\, X$ be any polystable vector bundle
with $c_1(W)\, =\, 0$. The Bogomolov inequality says that
\begin{equation}\label{bo}
\int_X c_2(W)\omega^{d-2}\, \geq\, 0
\end{equation}
(see \cite[Ch.~IV, \S~4, p. 114, Theorem (4.7)]{Ko}, \cite{Bo}).
In particular, for any polystable vector bundle $F$,
\begin{equation}\label{bo2}
\int_X c_2(End(F))\omega^{d-2}\, \geq\, 0
\end{equation}
(from \cite[p. 285, Theorem 3.18]{RR} it follows that the
vector bundle $End(F)\,=\, F\otimes F^*$ is also polystable).

{}From Eq. \eqref{S},
$$
c_2(End(\text{ad}(E_{G'})))\, =\, \sum_{i=1}^\ell
c_2(End(S_i/S_{i-1}))
$$
$$
 +\sum_{j=2}^{\ell} \sum_{i=1}^{j-1}
c_2(((S_i/S_{i-1})\otimes (S_j/S_{j-1})^*)\oplus
((S_i/S_{i-1})^*\otimes (S_j/S_{j-1})))\, .
$$
Note that $((S_i/S_{i-1})\otimes (S_j/S_{j-1})^*)\oplus
((S_i/S_{i-1})^*\otimes (S_j/S_{j-1}))$ is a polystable
vector bundle because both $S_i/S_{i-1}$ and $S_j/S_{j-1}$
are polystable with
$$
\mu (S_i/S_{i-1})\, =\, \mu (S_j/S_{j-1})
$$
(see \cite[p. 285, Theorem 3.18]{RR}). Also,
$$
c_1(((S_i/S_{i-1})
\otimes (S_j/S_{j-1})^*)\oplus ((S_i/S_{i-1})^*
\otimes (S_j/S_{j-1})))\, =\, 0
$$
because $((S_i/S_{i-1})
\otimes (S_j/S_{j-1})^*)\oplus ((S_i/S_{i-1})^*
\otimes (S_j/S_{j-1}))$ is self--dual.
Hence from Eq. \eqref{bo} and Eq. \eqref{bo2} we have
\begin{equation}\label{bo3}
\int_X c_2(End(\text{ad}(E_{G'})))\omega^{d-2}\, \geq\,
\int_X c_2(End(S_i/S_{i-1}))\omega^{d-2}
\end{equation}
for all $i\, \in \, [1\, ,\ell]$.

We have $c_1({\rm ad}(E_{G'}))\, =\, 0$ (see Remark \ref{adj.}).
Hence from Eq. \eqref{a} it follows immediately that
$$
\int_X c_2(End(\text{ad}(E_{G'})))\omega^{d-2}\, =\, 0\, .
$$
Therefore, from Eq. \eqref{bo3} and Eq. \eqref{bo2},
\begin{equation}\label{bo4}
\int_X c_2(End({\rm ad}(E_{L(Q)})))\omega^{d-2}
\,=\, \int_X
c_2(End(S_{(\ell+1)/2}/S_{(\ell-1)/2}))\omega^{d-2}\, =\, 0\, .
\end{equation}

We have $c_1({\rm ad}(E_{L(Q)}))\, =\, 0$
(see Remark \ref{adj.}). Therefore, from
Eq. \eqref{bo4} it follows immediately that Eq. \eqref{c3} holds.

Next we will show that for every character $\chi$ of $L(Q)$,
the associated line bundle $E_{L(Q)}(\chi)
\,\longrightarrow\, X$ is of degree zero.

Consider the adjoint action of $Q$ on the Lie
algebra ${\mathfrak g}'$ of $G'$. Let
$$
0\, =\, W'_0\, \subset\, W'_1 \, \subset\, \cdots \, \subset\,
W'_{n-1} \, \subset\, W'_n\, =:\, {\mathfrak g}'
$$
be the filtration of $Q$--modules constructed as in Lemma
\ref{lem1}. Therefore, the subgroup $R_u(Q)$ in Eq. \eqref{Q}
acts trivially on the direct sum
\begin{equation}\label{hw}
\widehat{W}\, :=\, \bigoplus_{i=1}^n W'_i/W'_{i-1}\, .
\end{equation}
Consequently, the action of $Q$ on $\widehat{W}$ factors
through an action of $L(Q)$ on $\widehat{W}$. Since $G'$ acts
faithfully on ${\mathfrak g}'$ (the center of $G'$ is
trivial), the action of $L(Q)$ on 
$\widehat{W}$ is also faithful.

Let
\begin{equation}\label{bw}
{\mathcal W}\, :=\, E_{L(Q)}(\widehat{W})\,\longrightarrow\, X
\end{equation}
be the vector bundle associated to the principal $L(Q)$--bundle
$E_{L(Q)}$ for the $L(Q)$--module $\widehat{W}$
in Eq. \eqref{hw}. From the construction of the reduction $E_Q$
it follows that the vector
bundle ${\mathcal W}$ is isomorphic to the direct sum
\begin{equation}\label{ds}
\bigoplus_{i=1}^\ell S_i/S_{i-1}
\end{equation}
associated to the filtration in Eq. \eqref{S}. We recall that for
each $i\, \in\, [1\, ,\ell]$, the quotient $S_i/S_{i-1}$ is a
polystable vector bundle of degree zero. Hence the vector
bundle ${\mathcal W}$, which is isomorphic to the vector bundle
in Eq. \eqref{ds}, is also polystable of degree zero. Further
note that as elements of the Grothendieck $K$--group, we have
$$
[{\mathcal W}]\, =\,\sum_{i=1}^\ell [S_i/S_{i-1}]\, =\,
[\text{ad}(E_{G'})]\, \in\, K(X)\, .
$$
In particular,
\begin{equation}\label{dc2}
c_1({\mathcal W})\,=\, c_1(\text{ad}(E_{G'}))\, =\, 0
\end{equation}
and
\begin{equation}\label{d2}
c_2({\mathcal W})\,=\, c_2(\text{ad}(E_{G'}))\, .
\end{equation}

Take any one--dimensional $L(Q)$--module $\mathcal M$.
We noted earlier that $L(Q)$ acts faithfully on the
$L(Q)$--module $\widehat{W}$ in Eq. \eqref{hw}.
Consequently, there are nonnegative integers $a$ and $b$
such that $\mathcal M$ is a direct summand of the $L(Q)$--module
\begin{equation}\label{dc3}
W_{a,b}\, := \, (\otimes^a \widehat{W})
\otimes (\otimes^b \widehat{W}^*)
\end{equation}
\cite[p. 40, Proposition 3.1(a)]{De}. Let
\begin{equation}\label{dc4}
{\mathcal W}_{a,b}\, :=\, E_{L(Q)}(W_{a,b})\,\longrightarrow\, X
\end{equation}
be the vector bundle associated to the principal $L(Q)$--bundle
$E_{L(Q)}$ for the $L(Q)$--module $W_{a,b}$ in Eq. \eqref{dc3}.
Let
\begin{equation}\label{dc5}
L_{\mathcal M}\, :=\, E_{L(Q)}({\mathcal M}) \,\longrightarrow\, X
\end{equation}
be the holomorphic
line bundle associated to the principal $L(Q)$--bundle
$E_{L(Q)}$ for the above one--dimensional
$L(Q)$--module $\mathcal M$. Since $\mathcal M$
is a direct summand of the $L(Q)$--module $W_{a,b}$, there
is a holomorphic vector bundle $V$ over $X$ such that
\begin{equation}\label{dc6}
{\mathcal W}_{a,b}\, =\, L_{\mathcal M}\oplus V\, ,
\end{equation}
where ${\mathcal W}_{a,b}$ is constructed in Eq. \eqref{dc4}.

{}From Eq. \eqref{dc3} and Eq. \eqref{bw} we have
\begin{equation}\label{dc7}
{\mathcal W}_{a,b}\, =\, (\otimes^a \mathcal{W})
\otimes (\otimes^b \mathcal{W}^*)\, .
\end{equation}
The vector bundle $\mathcal{W}^*$ is polystable of degree zero
because $\mathcal{W}$ is so. Hence from Eq.
\eqref{dc7} it follows that ${\mathcal W}_{a,b}$ is also
polystable of degree zero \cite[p. 285, Theorem 3.18]{RR}.
Now from Eq. \eqref{dc6} it follows immediately that
$$
\text{degree}(L_{\mathcal M})\, =\, 0\, .
$$

To complete the proof of the proposition we need to show that
Eq. \eqref{c03} holds.

{}From Eq. \eqref{dc2} and Eq. \eqref{dc7}
we have
\begin{equation}\label{dc9}
c_1({\mathcal W}_{a,b})\, =\, 0\, .
\end{equation}
In view of Eq. \eqref{dc2} and Eq. \eqref{dc7}, from Eq.
\eqref{d2} and Eq. \eqref{a} it follows that
\begin{equation}\label{dc8}
\int_X c_2({\mathcal W}_{a,b})\omega^{d-2}\, =\, 0\, .
\end{equation}

Let $\nabla$ be the Hermitian--Einstein connection on
the polystable vector bundle ${\mathcal W}_{a,b}$. From
Eq. \eqref{dc9} and Eq. \eqref{dc8} it follows that
$\nabla$ is flat (see \cite[Ch.~IV, \S~4, p. 114,
Theorem (4.7)]{Ko}). Therefore, the $(1\, ,0)$--part
$\nabla^{1,0}$ of $\nabla$ is a holomorphic connection
on the vector bundle ${\mathcal W}_{a,b}$;
see \cite{At} for holomorphic connections. Fix a
decomposition of ${\mathcal W}_{a,b}$ as in Eq.
\eqref{dc6}. Let
$$
\iota\, :\, L_{\mathcal M}\, \longrightarrow\, {\mathcal W}_{a,b}
$$
be the inclusion, and let $q_L\, :\, {\mathcal W}_{a,b}\, 
\longrightarrow\, L_{\mathcal M}$ be the projection given by this
decomposition. Now note that the composition
$$
L_{\mathcal M} \, \stackrel{\iota}{\longrightarrow}\,
{\mathcal W}_{a,b}\, \stackrel{\nabla^{1,0}}{\longrightarrow}\, 
{\mathcal W}_{a,b}\otimes\Omega^1_X\,
\stackrel{q_L\otimes\text{Id}_{\Omega^1_X}}{\longrightarrow}
\, L_{\mathcal M} \otimes\Omega^1_X\, ,
$$
where $\Omega^1_X$ is the holomorphic cotangent bundle
of $X$, is a holomorphic connection on the line 
bundle $L_{\mathcal M}$. Since $L_{\mathcal M}$ admits
a holomorphic connection, we have
$$
c_1(L_{\mathcal M})\, =\, 0
$$
(see \cite[p. 196,Proposition 12]{At}). In particular,
Eq. \eqref{c03} holds. This completes the proof of the
proposition.
\end{proof}

\section{Pseudostability and flatness}

We put down the previous results in the form of the following
theorem.

\begin{theorem}\label{thm1}
Let $E_G$ be a holomorphic principal $G$--bundle over a
compact connected K\"ahler manifold $(X\, , \omega)$, where $G$
is a connected reductive linear algebraic group
defined over $\mathbb C$. 
The following three statements are equivalent:
\begin{enumerate}

\item{} There is a parabolic subgroup $P\, \subset\, G$
and a holomorphic reduction of structure group
$E_P\, \subset\, E_G$ to $P$, such that the corresponding
$L(P)/Z(G)$--bundle
$$
E_P(L(P)/Z(G))\,=\, (E_P\times (L(P)/Z(G)))/P\,
\longrightarrow\, X
$$
admits a unitary flat connection.

\item{} The adjoint
vector bundle ${\rm ad}(E_G)$ is numerically flat.

\item{} The principal $G$--bundle $E_G$ is pseudostable,
and
$$
\int_X c_2({\rm ad}(E_G))\omega^{d-2}\, =\, 0\, .
$$
\end{enumerate}
\end{theorem}

\begin{proof}
Assume that the first statement holds.
As in Eq. \eqref{ep0}, set $L'(P)\, :=\, L(P)/Z(G)$. So
$E_{L'(P)}\, :=\, E_{L(P)}/Z(G)$ satisfies all the conditions
in Proposition \ref{prop1}. Hence from Corollary \ref{cor1}
we know that ${\rm ad}(E_G)$ is numerically flat.

{}From Lemma \ref{lem0} we know that the second statement 
implies the third statement.

Finally, assume that the third statement holds. Define
$G'\, =\, G/Z(G)$ and $E_{G'}$ as in Eq. \eqref{e0p} and Eq.
\eqref{e00} respectively. Using the decomposition in Eq.
\eqref{d} we conclude that
$$
\int_X c_2({\rm ad}(E_{G'}))\omega^{d-2}\, =\, 0\, .
$$
Therefore, from Proposition \ref{prop2} we know that
$E_L(Q)$ (constructed in Eq. \eqref{EL}) is
polystable, and the two statements in Proposition \ref{prop2}
hold.

Let $P$ be the inverse image of $Q\, \subset\, G'$ under the
quotient map $G\, \longrightarrow\, G'$.
We note that $P$ is a parabolic
subgroup of $G$. The inverse image of $E_Q\, \subset\, E_{G'}$ 
under the quotient map
$$
E_G\, \longrightarrow\, E_{G'}\,=\, E_G/Z(G)
$$
(see Eq. \eqref{e00}) is a holomorphic reduction of structure
group of $E_G$ to $P$. In the proof of Proposition
\ref{prop1} we saw that the Hermitian--Einstein connection on
$E_{L'(P)}$ is flat. Therefore, the first statement
in the theorem holds. This completes the proof of the theorem.
\end{proof}

Assume that $X$ is a complex projective manifold, and $\omega$
represents a rational cohomology class. We will show that the
third statement in Theorem \ref{thm1} is equivalent to the
first statement in Theorem \ref{thm-gb}.

Let $E_G$ be a holomorphic principal $G$--bundle over the
complex projective manifold $X$.
In \cite[p. 26, Proposition 2.4]{BG} it was shown that
if $c_2(\text{ad}(E_G))\, \in\, H^4(X,\, {\mathbb Q})$
vanishes, then $E_G$ is semistable if and only if
$E_G$ is pseudostable. Therefore, the
first statement in Theorem \ref{thm-gb} implies the
third statement in Theorem \ref{thm1}.

We now note that in the proof of
\cite[Proposition 2.4]{BG}, the weaker condition
that
$$
\int_X c_2({\rm ad}(E_G))\omega^{d-2}\, =\, 0
$$
is needed; see the last five lines in page 28 of \cite{BG}.
Note that Theorem 2 in \cite[p. 39]{Si} assumes only
$$
\int_X c_2({\rm ad}(E_G))\omega^{d-2}\, =\, 0
$$
and not the a priori stronger condition that
$c_2(\text{ad}(E_G))\,=\,0$.

Since $c_1({\rm ad}(E_G))\, =\, 0$ (see Remark
\ref{adj.}), from \cite[p. 39, Theorem 2]{Si} we
conclude that $c_2({\rm ad}(E_G))\, =\, 0$
if ${\rm ad}(E_G)$ is semistable with
$$
\int_X c_2({\rm ad}(E_G))\omega^{d-2}\, =\, 0\, .
$$

Therefore, we have the following corollary:

\begin{corollary}\label{cor2}
Assume that $X$ is a complex projective manifold, and the
cohomology class in $H^2(X,\, {\mathbb R})$ represented
by the K\"ahler form $\omega$ lies in $H^2(X,\, {\mathbb Q})$.
Then the third statement in Theorem \ref{thm1} is
equivalent to the first statement in Theorem \ref{thm-gb}.
\end{corollary}


\end{document}